\documentclass[11pt,
]{article}

\usepackage{amsmath}

\usepackage{amssymb}
\usepackage{latexsym}

\usepackage{graphicx,color}
\usepackage{times}

\usepackage[colorlinks=true, breaklinks=true, pdfstartview=FitV, linkcolor=red, citecolor=blue, urlcolor=black]
{hyperref}

\usepackage[T1]{fontenc}

\allowdisplaybreaks

\begin{document}

\title{Note on the analytical integration of circumterrestrial orbits\thanks{%
Submitted to Advances in Space Research}}

\author{Martin Lara%
\thanks{GRUCACI, University of La Rioja, and Space Dynamics Group -- UPM}
\thanks{\tt{mlara0@gmail.com} }
}

\date{\today }

\maketitle

\begin{abstract}
The acclaimed merits of analytical solutions based on a fictitious time developed in the 1970's were partially overvalued due to a common misuse of classical analytical solutions based on the physical time that were taken as reference. With the main problem of the artificial satellite theory as a model, we carry out a more objective comparison of both kinds of theories. We find that the proper initialization of classical solutions notably balances the performance of the two distinct approaches in what respects to accuracy. Besides, extension of both kinds of satellite theories to higher orders show additional pros and cons of each different perturbation approach, thus providing complementary information to prospective users on which kind of analytical solution may better support their needs.
\end{abstract}

\section{Introduction} \label{Introd}


The advent of the perturbation theories of Kozai \cite{Kozai1959} and Brouwer \cite{Brouwer1959} made a breakthrough in the analytical computation of circumterrestrial orbits. Both theories were computed in closed from of the eccentricity, thus extending their range of application from the lower eccentricities to highly elliptic orbits, and successfully included the effects of higher order zonal harmonics of the geopotential as well as the second order effects of the dominant, zonal harmonic of the second degree ---customarily denoted $J_2$. Brouwer and Kozai analytical solutions are given in the form of collections of truncated power series of $J_2$ with time as the argument. Several of the series represent the secular terms of the solution, which are expected to provide the average orbit evolution, and the remaining series comprise the periodic corrections that are needed for computing ephemeris. In the original proposals of Kozai and Brouwer, the former were accurate up to the second order of $J_2$ effects, whereas the later were truncated to $\mathcal{O}(J_2)$.
\par

A crucial point in the use of these kinds of truncated power series is the initialization of the constants of the theory \cite{Cain1962,Walter1967,Ustinov1967,EcksteinHechler1970,GaiasColomboLara2020}. In particular, the accurate determination of the ``mean'' mean motion is mandatory to obtain the in-track error that would be  expected from the truncation of the secular terms. This point was clearly highlighted by developers of perturbation solutions (see remarks in \cite{Kozai1962}). However, it created some confusion among users, who reported unexpected along-track errors in the propagation of initial conditions when Brouwer's theory was compared with alternative available solutions (see \cite{BonavitoWatsonWalden1969} and references therein). The trouble resulted from Brouwer's computations of the periodic corrections ---which are needed for the initialization of the constants of the theory from given initial conditions--- to a lower order than the one of the secular terms. On the other hand, this inaccuracy did not happen when the value of the ``mean'' mean motion was known to the required accuracy by other means, like the usual fit to observations. Therefore, the problem was understood like an inconsistency in the use of Brouwer's solution.
\par

Increasing the truncation order of the periodic corrections to solve the inconsistency involves heavy additional analytical computations \cite{Kozai1962,DepritRom1970}, which commonly discouraged users from their implementation. However, since the increased accuracy was fundamentally needed for the osculating to mean conversion of the semimajor axis, simple ways of patching Brouwer's solution were soon devised \cite{LyddaneCohen1962}. In this regard, the alternative proposed by Breakwell and Vagners \cite{BreakwellVagners1970} is remarkable, because it does not require additional computations to those already carried out by Brouwer. Indeed, the smart use of the energy equation to calibrate the mean semimajor axis provided a computationally inexpensive way of reducing in-track errors to a comparable level of those in the radial and cross-track directions (see \cite{Lara2020arxiv}, for instance).
\par

Conversely, because the accurate initialization of the mean semimajor axis is clearly associated with the correct treatment of the total energy, the pronounced loss of accuracy in the in-track direction does not happen to perturbation solutions formulated in the extended phase space, in which the total energy is one more integration variable \cite{Scheifele1970,Scheifele1970CeMec,StiefelScheifele1971}. In this scenario, the analytical solution is advantageously computed in a fictitious-time scale, whereas the physical time is incorporated into the perturbation scheme like the conjugate coordinate of the total energy. In particular, the annoying complications introduced by the equation of the center when dealing with geopotential perturbations in closed form \cite{Jefferys1971,Deprit1981,Metris1991,Ahmed1994,Lara2019CMDA} are now totally avoided.
\par

The merits of the extended phase space formulation in the implementation of an analytical orbit generator are unquestionable \cite{Scheifele1981}. However, assertions like the accuracy of first order analytical solutions based on this approach compare with traditional second order perturbation solutions based on the physical time \cite{ScheifeleGraf1974} are true only when the latter is incorrectly initialized, and hence must be somewhat downgraded. To show that we take the main problem of the artificial satellite theory as a demonstration model, and make an independent implementation of both kinds of theories. Needless to say that our implementation of the traditional one, based on the physical time, incorporates the calibration of the mean semimajor axis in the initialization part of the orbit theory. While we still found higher accuracy of the secular terms in the case of the fictitious-time approach, in no way it reaches the accuracy of a second-order traditional (correctly initialized) perturbation solution. In addition, we show that, due to the unavoidable truncation of the time equation, additional errors arise in the fictitious-time approach that are of comparable magnitude to the obtained position errors. Moreover, the need of inverting the time equation, or resorting to iterative procedures for finding the fictitious time corresponding to a given physical time \cite{Bond1979}, which is needed for usual ephemeris generation, must also be taken as a shortcoming of the fictitious-time approach because it notably increases the computational burden of usual ephemeris evaluation.
\par

To further assess the relative merits of each perturbation approach, we extended both kinds of satellite theories to the second order that is customarily used in accurate analytic ephemeris prediction programs \cite{CoffeyAlfriend1984,CoffeyNealSegermanTravisano1995}. In the case of the extended phase space formulation, spurious higher-order terms are generated at each order of the perturbation algorithm, a fact that makes the series defining the analytical solution notably longer than those resulting from the traditional, physical-time approach. Therefore, a post-processing is mandatory to keep the length of the series comprising the perturbation solution based on the fictitious time of comparable size to those of the traditional solution based on the physical time.
\par

\section{Solution to the $J_2$-problem in the extended phase space}

The topic of perturbation solutions of the artificial satellite problem using the traditional approach is satisfactorily covered in the literature. In particular, for the main problem, explicit expressions of higher-order solutions have been provided in different places (see \cite{Healy2000,Lara2019UR,Lara2020}, for instance). Therefore, we focus on the less studied formulation of perturbation solutions in the extended phase space, which, as customary nowadays, we recompute using the method of Lie transforms \cite{Deprit1969} rather than the von Zeipel perturbation algorithm used in the original developments \cite{ScheifeleGraf1974}.
\par

We approach the perturbation solution in the true anomaly Delaunay-similar elements \cite{Scheifele1970,Scheifele1970CeMec,Floria1997}. These elements are the action-angle variables of the unperturbed problem ---the Kepler problem in the extended phase space, with the true anomaly as the independent variable--- in which the complete Hamiltonian reduction of perturbed problems is naturally achieved. The fact that the Delaunay-similar elements are harmed by singularities is not at all of concern. Indeed, once the perturbation solution is computed, we can always reformulate it at our will in our preferred choice of singular or non-singular, canonical or non-canonical variables \cite{Lyddane1963,DepritRom1970}.
\par

\subsection{The main problem in Delaunay-similar variables}

The main problem admits Hamiltonian formulation. Using polar variables, the Hamiltonian is written in the form (see \cite{DepritFerrer1987} for instance)
\begin{equation} \label{HamP}
\mathcal{H}\equiv\frac{1}{2}\left(R^2+\frac{\Theta^2}{r^2}\right)-\frac{\mu}{r}+J_2\frac{\mu}{r}\frac{R_\oplus^2}{r^2}P_2(s\sin\theta),
\end{equation}
where $\mu$, $R_\oplus$, and $J_2$, are physical parameters that characterize the gravity field, and stand for the Earth's gravitational parameter, equatorial radius, and oblateness coefficient, respectively; $r$ is the radius, $R$ the radial velocity, the polar angle $\theta$ is commonly dubbed as the argument of the latitude, $\Theta$ is the total angular momentum (per mass unit), $P_2$ is the Legendre polynomial of degree 2, and $s\equiv\sin{I}$, with $I$ denoting orbital inclination. The third component of the specific angular momentum $N=\Theta\cos{I}$ is an integral that stems from the axial symmetry of the main problem.
\par

Hamiltonian (\ref{HamP}) is conveniently formulated in the \emph{extended phase space} \cite{Poincare1892vII,StiefelScheifele1971} by adding the pair of variables $q=t$ and $Q=-\mathcal{H}$, and defining the new Hamiltonian $\mathcal{H}^*\equiv\mathcal{H}+Q$, which, therefore, is constrained to the manifold $\mathcal{H}^*=0$. If, besides, we define a new time scale $\tau$, such that $\mathrm{d}t=\Psi(r,\theta,R,\Theta,q,Q)\,\mathrm{d}\tau$, it can be checked that the Hamiltonian $\widetilde{\mathcal{H}}\equiv\mathcal{H}^*\Psi$ remains constant in the new time scale, and that the solution of the system governed by $\widetilde{\mathcal{H}}$ coincides with the system governed by $\mathcal{H}^*$ when the initial conditions for $\tau=0$ are the same as the ones for $t=0$ \cite{Scheifele1970CeMec}. In particular, the choice
\begin{equation} \label{dttau}
\mathrm{d}t=\frac{r^2}{\Gamma(r,\theta,\Theta,N)}\,\mathrm{d}\tau,
\end{equation}
allows us to replace Eq.~(\ref{HamP}) by
\begin{equation} \label{Hampolar}
\widetilde{\mathcal{H}}\equiv\left[\frac{1}{2}\left(R^2+\frac{\Theta^2}{r^2}\right)-\frac{\mu}{r}+J_2\frac{\mu}{r}\frac{R_\oplus^2}{r^2}P_2(s\sin\theta)+Q\right]\frac{r^2}{\Gamma},
\end{equation}
the solutions of which only make sense in the manifold $\widetilde{\mathcal{H}}=0$.
\par

When choosing
\begin{equation} \label{GammaP}
\Gamma=\frac{1}{2}\Theta\left\{1+\left[1+2J_2\frac{\mu}{r}\frac{R_\oplus^2}{\Theta^2}P_2(s\sin\theta)\right]^{1/2}\right\},
\end{equation}
which only differs from $\Theta$ on the order of $J_2$, then Hamiltonian (\ref{Hampolar}) is suitably formulated in the Delaunay-similar variables $(\phi,g,h,\lambda,\Phi,G,H,\Lambda)$, where $\phi$ is the true anomaly, $g$ is the argument of the perigee, $h$ is the right ascension of the ascending node, $\lambda$ is the time element, $\Phi$ is related to the Keplerian energy, $G=\Theta$, $H=N$, and $\Lambda=Q$. They are computed from the polar variables following the sequence (see \cite{ScheifeleGraf1974,Floria1997} for details)
\begin{align}
\Phi=& \; 2(\Theta-\Gamma)+\frac{\mu}{(2Q)^{1/2}}, \\
p=& \; \frac{1}{\mu}(2\Theta-\Gamma)^2, \\
e=& \; \sqrt{1-2Qp/\mu}, \\
\cos\phi=& \; \frac{1}{e}\left(-1+\frac{p}{r}\right),\qquad \sin\phi=\frac{R}{e}\sqrt{\frac{p}{\mu}}, \\
u=& \; 2\arctan\left(\sqrt{\frac{1-e}{1+e}}\tan\frac{\phi}{2}\right), \\ \label{lambda}
\lambda=& \; q-\frac{\mu}{(2\Lambda)^{3/2}}\left(u-e\sin{u}-\phi\right).
\end{align}
In these variables, Eq.~(\ref{Hampolar}) takes the form
\begin{equation} \label{Famain}
\mathcal{F}\equiv \Phi-\frac{\mu}{\sqrt{2\Lambda}}
-J_{2}\frac{\mu}{r}\frac{R_\oplus^2}{\Gamma}\frac{1}{4}\left[2-3s^2+3s^2\cos2(g+\phi)\right],
\end{equation}
in which, now, $s^2=1-H^2/G^2$, $r=p/(1+e\cos\phi)$, and
\begin{equation}
\Gamma=G-\frac{1}{2}\left(\Phi-\frac{\mu}{\sqrt{2\Lambda}}\right).
\end{equation}
\par

\subsection{Perturbation approach}

Due to the smallness of the Earth's $J_2$ coefficient, the complete Hamiltonian reduction of the main problem Hamiltonian (\ref{Famain}) can be achieved, up to some truncation order, by the Lie transforms method \cite{Deprit1969,BoccalettiPucacco1998v2,Lara2021}. We assume that the reader is enough familiarized with this perturbation method and only provide the results.
\par

In our approach, we adhere to the tradition and carry out the sequential elimination of short- and long-period terms \cite{Brouwer1959,ScheifeleGraf1974}. This is done analytically by changing first from original to prime variables, and then from prime variables to mean (double-prime) variables. The complete Hamiltonian reduction is then achieved by neglecting higher-order terms from the mean elements Hamiltonian, to obtain
\begin{equation} \label{Fprime}
\mathcal{F}''\equiv\mathcal{F}''(-,-,-,-,\Phi'',G'',H'',\Lambda'').
\end{equation}
The analytical solution is trivial in mean elements. Indeed, the Hamilton equations of Eq.~(\ref{Fprime}) immediately show that the momenta are constant in mean elements
\begin{equation} \label{momsec}
\Phi''=\Phi''_0, \quad
G''=G''_0, \quad
H''=H, \quad
\Lambda''=\Lambda,
\end{equation}
where $H$ and $\Lambda$ are integrals of the original main problem, whereas the angles evolve linearly with constant frequencies given by
\begin{align}
n_\phi= & \; n_\phi(\Phi''_0,G''_0,H,\Lambda)={\partial\mathcal{F}''}/{\partial\Phi''}, \\
n_g= & \; n_g(\Phi''_0,G''_0,H,\Lambda)={\partial\mathcal{F}''}/{\partial{G}''}, \\
n_h= & \; n_h(\Phi''_0,G''_0,H,\Lambda)={\partial\mathcal{F}''}/{\partial{H}''}, \\ \label{nlambda}
n_\lambda= & \; n_\lambda(\Phi''_0,G''_0,H,\Lambda)={\partial\mathcal{F}''}/{\partial\Lambda''}.
\end{align}
\par

To avoid offending divisors in the perturbation series, the analytical solution is customarily formulated in the canonical set of non-singular Poicar\'e-similar elements \cite{Scheifele1981}. However, due to the axial symmetry of the main problem model, we find convenience in formulating the solution in a non-canonical set of variables that replaces the troublesome elements by the argument of the latitude $\theta=\phi+g$ and the semi-equinoctial elements defining the eccentricity vector in the orbital plane $C=e\cos{g}$, $S=e\sin{g}$. The solution in mean elements is thus
\begin{align} \label{lambdasec}
\lambda''= & \; \lambda_0''+n_\lambda\tau, \\
\theta''= & \; \theta_0''+(n_\phi+n_g)\tau, \\
C''= & \; C''_0\cos(n_g\tau)-S''_0\sin(n_g\tau), \\
S''= & \; S''_0\cos(n_g\tau)+C''_0\sin(n_g\tau), \\  \label{hsec}
h''= & \; h''_0+n_h\tau,
\end{align}
and the original dynamics is recovered by plugging this solution into the transformation from mean to osculating elements, which in the current, Hamiltonian case is derived from a generating function.
\par

We extended the solution of \cite{ScheifeleGraf1974} to the second order. That is, the secular terms in Eqs.~(\ref{lambdasec})--(\ref{hsec}) are accurate to $\mathcal{O}(J_2^3)$ effects, whereas the periodic corrections are accurate only to $\mathcal{O}(J_2^2)$. Accordingly, the analytical solution comprises 33 perturbation series in addition to those defining the secular frequencies. Namely, $5\times3$ for the long-period corrections (first-order, second-order direct, and second-order inverse), and $6\times3$ for the short-period corrections (first-order, second-order direct, and second-order inverse). Remarkably, there is no need of integrating Eq.~(\ref{dttau}) to recover the physical time, which, on the contrary, due to the extended phase space formulation is obtained by making explicit $q$ from Eq.~(\ref{lambda}).
\par

Note that, due to the truncation of the periodic corrections, the initialization of the perturbation theory from a given set of initial conditions provides the needed constants in Eqs.~(\ref{momsec})--(\ref{nlambda}) only up to $\mathcal{O}(J_2^2)$ effects, whereas the secular frequencies are expected to be accurate to the order of $J_3^3$. While this is an inconsistency in traditional perturbation theories, which need to be patched to avoid large in-track errors, we will see that this uneven truncation is not of concern in the extended phase space formulation.

\subsection{First-order solution}

The generating function of the short-period elimination is
\begin{align} \nonumber
\mathcal{W}_1= & \;-\frac{1}{8}\Gamma\frac{R_\oplus^2}{\rho^2}\Big[
(4-6s^2)e\sin\phi +3es^2\sin(2g+\phi)  \\ \label{W1}
& \; +3s^2\sin(2g+2\phi)+es^2\sin(2g+3\phi) \Big],
\end{align}
where we defined the auxiliary variable
\begin{equation} \label{rho}
\rho=\Gamma\sqrt{p/\mu}.
\end{equation}
Remarkably, Eq.~(\ref{W1}) is formally the same as the first order term of the generating function of the elimination of the parallax transformation when it is written in Delaunay variables (cf.~\cite{LaraSanJuanLopezOchoa2013b}).
\par

Up to the first order, the transformation of the short-period elimination is computed like $\xi=\xi'+J_2\delta_{1,\xi}$, with $\xi$ denoting any of the involved variables. The short-period corrections $\delta_{1,\xi}$ are obtained from the computation of the Poisson bracket $\{\xi,\mathcal{W}_1\}$, and reformulating the result in prime variables. The inverse transformation is $\xi'=\xi-J_2\delta_{1,\xi}$ where now $\delta_{1,\xi}$ is written in the original, non-primed variables.
\par

For the long-period elimination we obtain the first order generating function
\begin{equation} \label{V1}
\mathcal{V}_{1}=\Gamma\frac{R_\oplus^2}{\rho^2}\frac{3}{32}
\frac{1}{\Delta}\left[15s^2-14+12(s^2-1)\delta\right]s^2e^2\sin2g,
\end{equation}
which must be written in the prime variables. In Eq.~(\ref{V1}) we abbreviated
\begin{equation} \label{icrit}
\Delta=3(5s^2-4)+6(s^2-1)\delta+2(3s^2-2)\upsilon,
\end{equation}
and introduced the auxiliary, non-dimensional functions
\begin{equation} \label{hidden}
\delta=-1+\Gamma/G, \qquad \upsilon=-1+\rho/p,
\end{equation}
which are $\mathcal{O}(J_2)$, and hence produce a small displacement of similar magnitude in the critical inclination value, which now happens for $\Delta=0$. That is, $\sin{I}=\sqrt{4/5}+\mathcal{O}(J_2)$. The first-order transformation to mean (double-prime) elements is computed like before. Thus, $\xi'=\xi''+J_2\delta'_{1,\xi'}$, where the long-period corrections $\delta'_{1,\xi'}$ are obtained computing the Poisson bracket $\{\xi',\mathcal{V}_1\}$, and reformulating the result in the double-prime variables. Analogously, $\xi''=\xi'-J_2\delta'_{1,\xi'}$ with $\delta'_{1,\xi'}$ now written in prime variables.
\par

The secular frequencies of the first order theory are obtained from the Hamilton equations of the secular Hamiltonian
\begin{equation} \label{Fsecond}
\mathcal{F}''=\Phi-\frac{\mu}{\sqrt{2\Lambda}}+\sum_{m\ge1}\frac{J_2^m}{m!}\mathcal{F}_m,
\end{equation}
where
\begin{equation} \label{F01}
\mathcal{F}_{1}=\Gamma\frac{R_\oplus^2}{\rho^2}\frac{1}{4}(3s^2-2),
\end{equation}
and
\begin{align} \nonumber
\mathcal{F}_{2}= & \; 
\Gamma\frac{R_\oplus^4}{\rho^4}\frac{1}{64}\big[ 4(15s^4-6s^2-4)+3(5s^4+8s^2-8)e^2
+24s^2  \\ \label{Q02}
& \; \times(2e^2+3)(s^2-1)\delta -2(e^2+1)(15s^4-24s^2+8)\upsilon \big],
\end{align}
are written in mean (double-prime) variables.

\subsection{Second-order solution}

The second order generating function takes the form
\begin{align} \nonumber
\mathcal{W}_2= & \;
\Gamma\frac{R_\oplus^4}{\rho^4}\frac{1}{3840}\Big\{
 60\big[45s^4+72s^2-80+168s^2(s^2-1)\delta-4  \\ \nonumber
& \;  \times(33s^4-48s^2+16)\upsilon\big]e\sin\phi 
+360\big[(5s^2-4)s^2+4s^2  \\ \nonumber
& \; \times(s^2-1)\delta\big]e^2\sin2\phi
-90\big[225s^2-206 +168(s^2-1)\delta \\ \nonumber
& \; +4(3s^2-2)\upsilon\big]s^2e\sin(2g+\phi)
-120\big\{39s^2-38+36(s^2  \\ \nonumber
& \;  -1)\delta-2(3s^2-2)\upsilon +2e^2\big[3s^2-4+6(s^2-1)\delta-(3s^2  \\ \nonumber
& \;   -2)\upsilon\big] \big\}s^2\sin(2g+2\phi)
+10\big[75s^2-42-24(s^2-1)\delta \\ \nonumber
& \; +28(3s^2-2)\upsilon\big] s^2e\sin(2g+3\phi)
+ 30\big[15s^2-14+12(s^2 \\ \nonumber
& \; -1)\delta\big]s^2e^2\sin(2g+4\phi) +45(4\upsilon+5)s^4e\sin(4g+3\phi)  \\ \nonumber
& \; + 45\big[2(\upsilon+1)+(2\upsilon+3)e^2\big]s^4\sin(4g+4\phi) \\
& \; +9 (4\upsilon+5)s^4e\sin(4g+5\phi)
\Big\}.
\end{align}
Now
\begin{equation} \label{short2nd}
\xi=\xi'+J_2\delta_{1,\xi}+\frac{1}{2}J_2^2\delta_{2,\xi},
\end{equation}
where $\delta_{2,\xi}$ is obtained after reformulating in prime variables the result of $\big\{\{\xi,\mathcal{W}_1\},\mathcal{W}_1\big\}+\{\xi,\mathcal{W}_2\}$. For the inverse transformation, 
\begin{equation} \label{short2ndi}
\xi'=\xi-J_2\delta'_{1,\xi}+\frac{1}{2}J_2^2\delta'_{2,\xi},
\end{equation}
where $\delta'_{2,\xi}$ is no longer the opposite of $\delta_{2,\xi}$, but the result, in the original variables, of $\big\{\{\xi,\mathcal{W}_1\},\mathcal{W}_1\big\}-\{\xi,\mathcal{W}_2\}$.
\par

The secular frequencies of the second-order theory are computed from the Hamiltonian (\ref{Fsecond}), which is now complemented with the third-order term
\begin{equation} \label{F3mean}
\mathcal{F}_{3}=\Gamma\frac{R_\oplus^6}{\rho^6}\frac{3}{2^{10}}
\frac{1}{\Delta^2}
\sum_{k=0}^{2}\sum_{j=0}^4\sum_{i=0}^{4-j}q_{k,i,j}(1+\delta)^i(1+\upsilon)^je^{2k},
\end{equation}
where the coefficients $q_{k,i,j}$ are given in Table~\ref{t:F3mean}.
\par

\begin{table*}[htbp]
\small \centering \tabcolsep 3.5 pt
\centerline{
\begin{tabular}{@{}llll@{}}
$i,j$ & \multicolumn{1}{c}{$k=0$} & \multicolumn{1}{c}{$k=1$} & \multicolumn{1}{c}{$k=2$} \\
\hline\noalign{\smallskip}
${}_{0,0}$ & $-72(2-3s^2)^3s^4$ & $-4(2-3s^2)^3 (41 s^4-24 s^2+8)$ & $18(2-3s^2)^3 s^4$ \\[0.3ex]
${}_{0,1}$ & $-4(2-3s^2)^3 (67 s^4-24 s^2+8)$ & $-2(2-3s^2)^3 (413 s^4-432 s^2+144)$ & $54(2-3s^2)^3 s^4$ \\[0.3ex]
${}_{0,2}$ & $-32(2-3s^2)^3 (15 s^4-24 s^2+8)$ & $-4(2-3s^2)^3 (513 s^4-816 s^2+272)$ & $9(2-3s^2)^3 s^4$ \\[0.3ex]
${}_{0,3}$ & $-16(2-3s^2)^3 (63 s^4-120 s^2+40)$ & $-16(2-3s^2)^3 (217 s^4-360 s^2+120)$ & $0$ \\[0.3ex]
${}_{0,4}$ & $-64(2-3s^2)^3 (17 s^4-24 s^2+8)$ & $-160(2-3s^2)^3 (17 s^4-24 s^2+8)$ & $0$ \\[0.3ex]
${}_{1,0}$ & $-144c^2 (2-3s^2)^2 s^2 (15 s^2-2)$ & $-64c^2 (2-3s^2)^2 (69 s^4-16 s^2+6)$ & $36c^2(2-3s^2)^2s^2(15 s^2-1)$ \\[0.3ex]
${}_{1,1}$ & $-48c^2(2-3s^2)^2(117s^4-24s^2+8)$ & $-12c^2(2-3s^2)^2(1135s^4-448s^2+224)$ & $72c^2(2-3s^2)^2s^2(18s^2-1)$ \\[0.3ex]
${}_{1,2}$ & $-288c^2 (2-3s^2)^2 (13 s^4-12 s^2+8)$ & $-96c^2 (2-3s^2)^2 (149 s^4-144 s^2+80)$ & $216c^2 (2-3s^2)^2 s^4$ \\[0.3ex]
${}_{1,3}$ & $-768c^2 (2-3s^2)^2 (3 s^4-6 s^2+4)$ & $-64c^2 (2-3s^2)^2 (141 s^4-224 s^2+120)$ & $0$ \\[0.3ex]
${}_{2,0}$ & $432c^2(2-3s^2)s^2(63s^4-76s^2+16)$ & $144c^2(2-3s^2)(399s^6-449s^4+80s^2-8)$ & $-27c^2(2-3s^2)s^2(243s^4-274s^2+40)$ \\[0.3ex]
${}_{2,1}$ & $144c^2(2-3s^2)(375s^6-411s^4+80s^2-8)$ & $144c^2(2-3s^2)(883s^6-947s^4+192s^2-40)$ & $-216c^2(2-3s^2)s^2(51s^4-56s^2+8)$ \\[0.3ex]
${}_{2,2}$ & $576c^2(2-3s^2)(36s^6-27s^4+8s^2-8)$ & $288c^2(2-3s^2)(223s^6-235s^4+96s^2-40)$ & $1296c^4(2-3s^2)s^4$ \\[0.3ex]
${}_{3,0}$ & $5184c^4s^2(5s^2-2)(6s^2-5)$ & $1728c^4s^2(205s^4-239s^2+56)$ & $-1296c^4s^2(30s^4-35s^2+8)$ \\[0.3ex]
${}_{3,1}$ & $3456c^4s^2(4s^2-3)(13s^2-4)$ & $1728c^4s^2(261s^4-289s^2+72)$ & $2592c^4(2-3s^2)s^2(5s^2-2)$ \\[0.3ex]
${}_{4,0}$ & $-15552c^6 s^2 (7 s^2-4)$ & $-38016c^6 s^2 (7 s^2-4)$ & $3888c^6 s^2 (7 s^2-4)$ \\
\noalign{\smallskip}\hline
\end{tabular}
}
\caption{Inclination polynomials $q_{k,i,j}$ in Eq.~(\protect\ref{F3mean}); $c=(1-s^2)^{1/2}$. }
\label{t:F3mean}
\end{table*}

The second-order term of the generating function of the transformation of the long-period elimination is
\begin{align} \nonumber
\mathcal{V}_{2}= & \; \Gamma\frac{R_\oplus^4}{\rho^4}\frac{3}{2^{10}}
\frac{1}{\Delta^3}
\sum_{l=1}^2\sum_{k=0}^{2-l}\sum_{j=0}^4\sum_{i=0}^{4-j}b_{l,k,i,j} \\ \label{V2}
& \; \times(1+\delta)^i(1+\upsilon)^je^{2k+2l}s^{2l}\sin2lg,
\end{align}
where the inclination polynomials $b_{l,k,i,j}$ are listed in Table~\ref{t:V2}. The transformation of the long-period elimination is obtained analogously to Eq.~(\ref{short2nd}) ---resp.~(\ref{short2ndi})--- and following formulas, replacing corresponding functions and variables by those of the long-period case.
\par

Remarkably, at this order there is an explosion in the number of terms of the transformation of the long-period elimination. This is clearly shown by comparison with Brouwer's 2nd order generating function of the long-period elimination, which takes the much simpler form
\begin{equation} \label{V2B}
\mathcal{V}^*_{2}=G\frac{R_\oplus^4}{p^4}\frac{1}{2^{10}}
\frac{1}{\tilde{\Delta}^2}\sum_{j=1}^2\frac{(1+\eta)^{j-2}}{\tilde{\Delta}^{j-1}}\sum_{i=0}^{3j^*}b_{j,i}\eta^je^{2j}s^{2j}\sin2jg,
\end{equation}
where $\tilde{\Delta}=5s^2-4$, $j^*=(j \bmod 2)$, and
\begin{equation}
\begin{array}{l}
b_{1,0}=-2(3975 s^6-6870 s^4+2928 s^2+16), \\
b_{1,1}= 2(1425 s^6-5370 s^4+6288 s^2-2320), \\
b_{1,2}= -2(15 s^2-14)(195 s^4-388 s^2+184), \\
b_{1,3}= 2(15 s^2-14)(45 s^4+36 s^2-56), \\
b_{2,0}= (15 s^2-14)^2(15 s^2-13).
\end{array}
\end{equation}
Still, if third-order corrections are not of concern, after computing the second-order corrections stemming from Eq.~(\ref{V2}), they can be notably simplified by making $\delta=\upsilon=0$. Other terms of the different series that compose the solution in the extended phase space are certainly amenable of analogous simplifications.
\par

\begin{table*}[htbp]
\small \centering \tabcolsep 3.5 pt
\centerline{
\begin{tabular}{@{}llll@{}}
$i,j$ & \multicolumn{1}{c}{$l=1,k=0$} & \multicolumn{1}{c}{$l=1,k=1$} & \multicolumn{1}{c}{$l=2,k=0$} \\
\hline\noalign{\smallskip}
${}_{0,0}$ & $8 (2-3s^2)^2(3 s^4+24 s^2-8)$ & $6 (2-3s^2)^2(3 s^4-24 s^2+8)$ & $-3 (2-3s^2)^3$ \\[0.3ex]
${}_{0,1}$ & $-24 (2-3s^2)^2(13 s^4-60 s^2+20)$ & $24 (2-3s^2)^2(7 s^4-24 s^2+8)$ & $-6 (2-3s^2)^3$ \\[0.3ex]
${}_{0,2}$ & $-32 (2-3s^2)^2(60 s^4-123 s^2+41)$ & $12 (2-3s^2)^2(37 s^4-72 s^2+24)$ & $-9 (2-3s^2)^3$ \\[0.3ex]
${}_{0,3}$ & $-8 (2-3s^2)^2(399 s^4-600 s^2+200)$ & $24 (2-3s^2)^2(15 s^4-24 s^2+8)$ & $0$ \\[0.3ex]
${}_{0,4}$ & $-48 (2-3s^2)^2(33 s^4-48 s^2+16)$ & $0$ & $0$ \\[0.3ex]
${}_{1,0}$ & $24 c^2 (2-3s^2) s^2(159 s^2-68)$ & $-12 c^2 (2-3s^2)(15 s^4+212 s^2-88)$ & $-90 c^2 (2-3s^2)^2$ \\[0.3ex]
${}_{1,1}$ & $24 c^2 (2-3s^2)(207 s^4+100 s^2-136)$ & $24 c^2 (2-3s^2)(81 s^4-384 s^2+152)$ & $-180 c^2 (2-3s^2)^2$ \\[0.3ex]
${}_{1,2}$ & $-72 c^2 (2-3s^2)(229 s^4-368 s^2+168)$ & $24 c^2 (2-3s^2)(279 s^4-560 s^2+200)$ & $-216 c^2 (2-3s^2)^2$ \\[0.3ex]
${}_{1,3}$ & $-192 c^2 (2-3s^2)(117 s^4-164 s^2+60)$ & $288 c^2 (2-3s^2)(15 s^4-24 s^2+8)$ & $0$ \\[0.3ex]
${}_{2,0}$ & $-144 c^2(552 s^6-1087 s^4+680 s^2-136)$ & $96 c^2(126 s^6-60 s^4-127 s^2+64)$ & $27 c^2 (2-3s^2)(35 s^2-34)$ \\[0.3ex]
${}_{2,1}$ & $-144 c^2(1065 s^6-2091 s^4+1324 s^2-256)$ & $192 c^2(45 s^6+87 s^4-217 s^2+88)$ & $-1728 c^4 (2-3s^2)$ \\[0.3ex]
${}_{2,2}$ & $-432 c^2 s^2(61 s^4-61 s^2+16)$ & $288 c^4(57 s^4-128 s^2+56)$ & $-1296 c^4 (2-3s^2)$ \\[0.3ex]
${}_{3,0}$ & $864 c^4(211 s^4-296 s^2+100)$ & $-1728 c^4 s^2(22 s^2-19)$ & $648 c^4(7 s^2-6)$ \\[0.3ex]
${}_{3,1}$ & $3456 c^4(62 s^4-85 s^2+32)$ & $-1152 c^4(27 s^4-17 s^2-4)$ & $-5184 c^6$ \\[0.3ex]
${}_{4,0}$ & $-10368 c^6(13 s^2-10)$ & $6912 c^6(5 s^2-2)$ & $-3888 c^6$ \\[0.3ex]
\noalign{\smallskip}\hline
\end{tabular}
}
\caption{Inclination polynomials $b_{l,k,i,j}$ in Eq.~(\protect\ref{V2}); $c=(1-s^2)^{1/2}$. }
\label{t:V2}
\end{table*}

\section{Accuracy tests}

A number of tests has been conducted to check the accuracy of the analytical solution based on the fictitious time and compare it with the precision provided by the traditional approach based on the physical time. For the later, the different analytical solutions in the literature for the same truncation order are essentially different arrangements of the same solution, which, therefore enjoy analogous accuracy \cite{Lara2020arxiv}. However, one formulation may be preferred over others depending on the user needs \cite{Lara2021IAC}. For our comparisons we choose the approach based on the sequential elimination of the parallax \cite{Deprit1981,LaraSanJuanLopezOchoa2013b}, the elimination of the perigee \cite{AlfriendCoffey1984,LaraSanJuanLopezOchoa2013c}, and the Delaunay normalization \cite{Deprit1982}, which seems to be a popular option in the literature \cite{CoffeyAlfriend1984,SanJuan1994,CoffeyNealSegermanTravisano1995,Lara2019UR}. Remark that, as opposite to common implementations by different authors, we always improve the initialization of the constants of the traditional perturbation solution with the calibration of the mean semimajor axis from the energy equation \cite{BreakwellVagners1970}.
\par

For the different orbital configurations tested, we always obtained similar comparative results between the traditional and fictitious-time modalities. Therefore, we only illustrate them for the PRISMA orbit \cite{PerssonJacobssonGill2005}, which is a low-eccentricity, sun-synchronous orbit at around 700 km altitude. In particular, in the computation of the initial conditions we used the orbital parameters $a=6878.14$ km, $e=0.001$, $I=97.42^\circ$, $\Omega=168.2^\circ$, $\omega=20^\circ$, and $M=30^\circ$. The true, reference orbit was then computed from the numerical integration of the main problem in Cartesian coordinates. To guarantee that all the stored digits of the reference orbit are exact in the floating-point number representation, the integration was carried out in extended precision. In addition, the integration of the fictitious time variation, given by Eq.~(\ref{dttau}), was incorporated to the differential system in order to check the accuracy provided by each theory exactly at the same physical time.
\par

\begin{figure}[htb]
\centering
\includegraphics[scale=0.7]{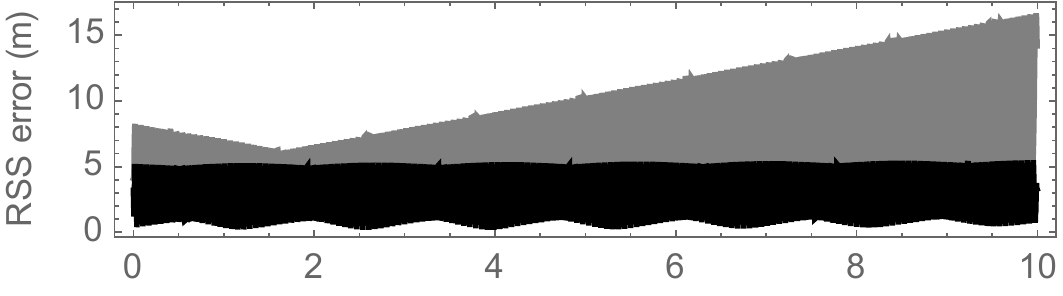}
\caption{RSS position errors of the first order theory in the extended phase space (black) superimposed to corresponding errors of the traditional approach (gray). Abscissas are days.}
\label{fig:RSS1st}
\end{figure}

Results showing the Root Sum Square (RSS) of the errors provided by the first-order truncation of each perturbation solution are shown superimposed in Fig.~\ref{fig:RSS1st} for a propagation interval of 10 days. The small secular trend in the errors of the traditional solution, of about 1 meter per day, is mostly nullified in the case of the solution based on the extended phase formulation, where the much smaller secular trend remains buried under the periodic oscillations of the errors. 
\par

\begin{figure*}[htb]
\centerline{
\includegraphics[scale=0.7]{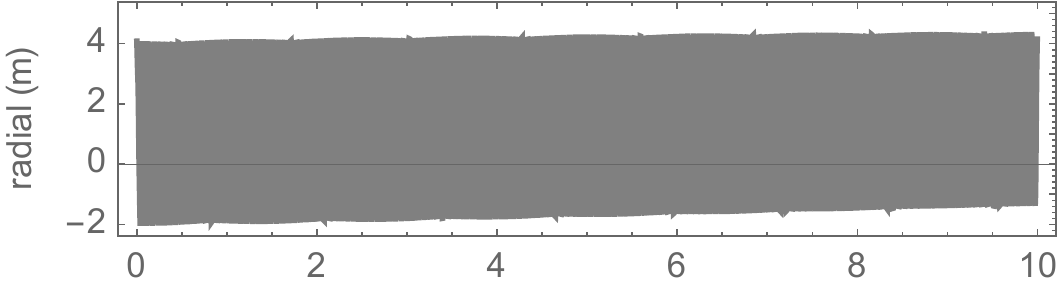}\qquad \includegraphics[scale=0.7]{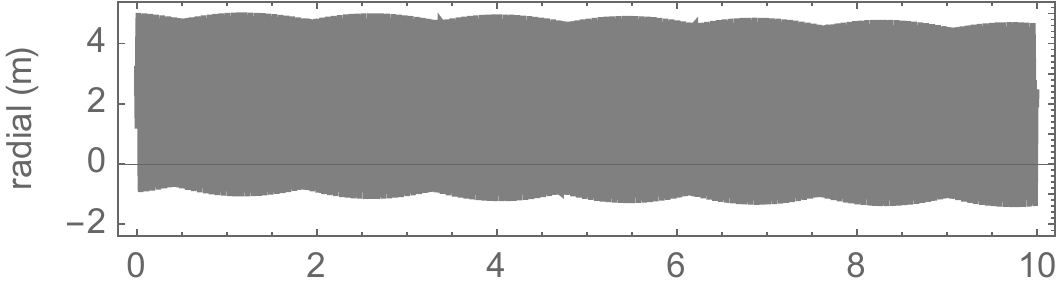} }
\centerline{
\includegraphics[scale=0.7]{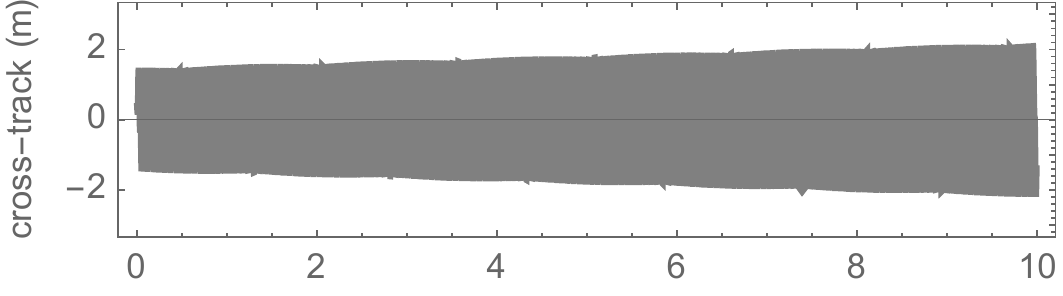} \qquad \includegraphics[scale=0.7]{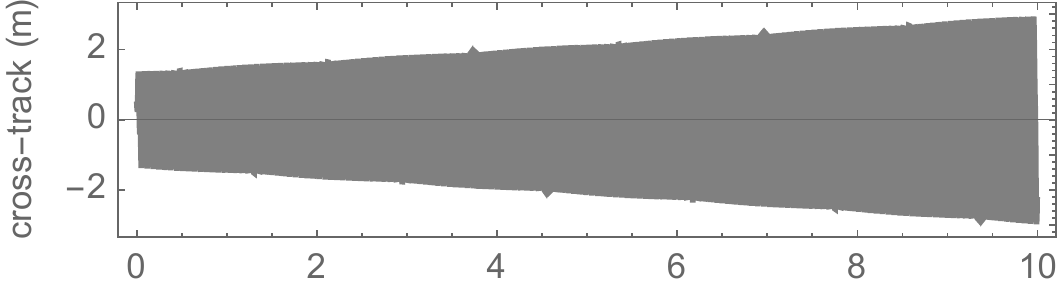} }
\centerline{
\includegraphics[scale=0.7]{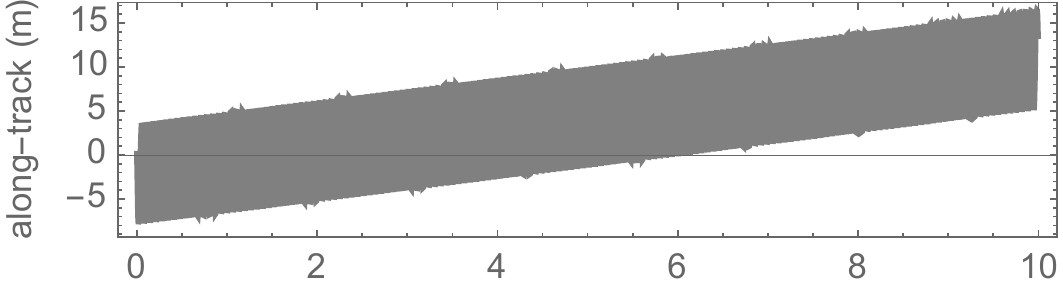} \qquad \includegraphics[scale=0.7]{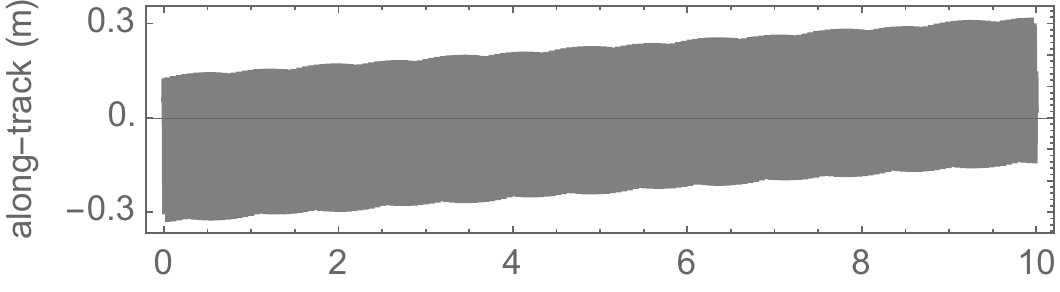} }
\caption{Intrinsic errors of the first-order traditional (left) and extended phase space theory (right). Abscissas are days.}
\label{fig:121intrinsic}
\end{figure*}

The nature of the position errors of each analytical solution is best illustrated when they are projected in the radial, along-track, and cross-track directions. Indeed, by simple inspection of Fig.~\ref{fig:121intrinsic}, we check that the errors in the radial and cross-track directions are periodic in nature and of comparable magnitude in both theories ---as it should be expected for perturbation solutions truncated to the same order. On the contrary, while the amplitude of the periodic components of the errors in the along-track direction remain in the case of traditional theory of the same level of accuracy as the other components of the errors, they reduce by about one order of magnitude, from meters to decimeters, in the case of the extended phase space formulation. The better performance of the extended phase space solution regarding along-track errors applies also to their evident secular component (bottom plots of Fig.~\ref{fig:121intrinsic}), which are mostly responsible from the different behavior of both theories, and clearly illustrate the benefits of including the total energy among the integration variables regarding the stabilization of the secular errors in the along-track direction.
\par

\begin{figure}[htb]
\centering
\includegraphics[scale=0.7]{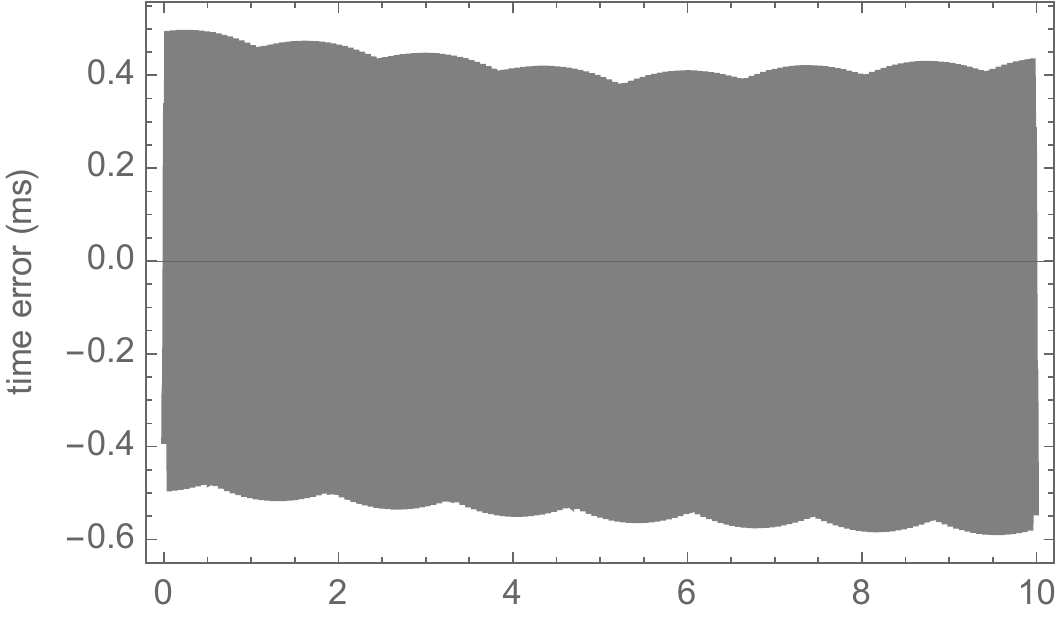}
\caption{Errors in the physical time determination stemming from the first-order theory in the extended phase space. Abscissas are days. }
\label{fig:timerror1}
\end{figure}

Undeniably, the greater secular trend of the errors of the traditional solution seems to make the extended phase space solution preferable in much longer propagation times. Still, for actual ephemeris computation, the fictitious time in which the extended phase space solution must be evaluated corresponding to the desired physical time will not be known in advance. On the contrary, its accurate computation requires a root-finding procedure that is commonly approached by Newton-Raphson iterations \cite{Bond1979}. Moreover, due to the unavoidable truncation of the perturbation solution, the fictitious time corresponding to a given physical time cannot be determined exactly. As shown in Fig.~\ref{fig:timerror1}, the nature of these errors is mostly periodic, with an amplitude of about one half of a millisecond. For the semimajor axis of the PRISMA orbit, this error in the timing may yield an additional in-track error in the meter level.
\par

\begin{figure}[htb]
\centering
\includegraphics[scale=0.7]{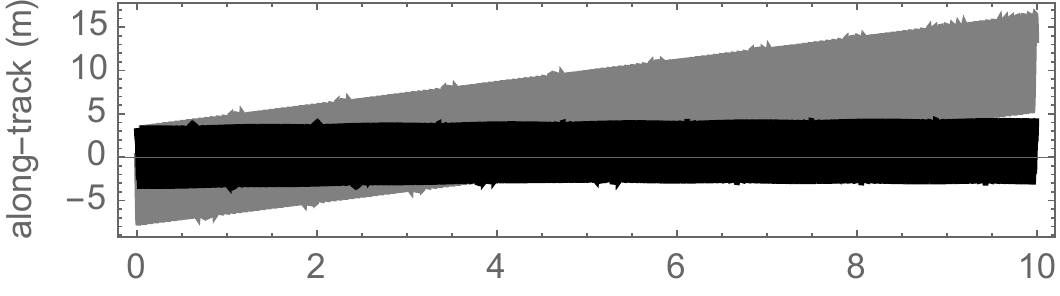}
\caption{Along-track errors of the extended phase space theory with the physical time as argument (black) superimposed to the bottom-left plot of Fig.~\protect\ref{fig:121intrinsic}. Abscissas are days. }
\label{fig:intrack121time}
\end{figure}

These additional inaccuracies derived from the $\mathcal{O}(J_2)$ truncation of the series that provide the physical time as a function of the fictitious one, which were not presented in Fig.~\ref{fig:RSS1st}, are now illustrated in Fig.~\ref{fig:intrack121time} where the along-track errors of the first-order theory in the extended phase space (in black) are computed taking the physical time as the argument, from which the needed fictitious time is then determined. This process requires the evaluation of the analytical theory in each iteration, thus notably increasing the computational burden, and, therefore, degrading the performance of the perturbation theory regarding computing time.
\par

\begin{figure}[htb]
\centering
\includegraphics[scale=0.7]{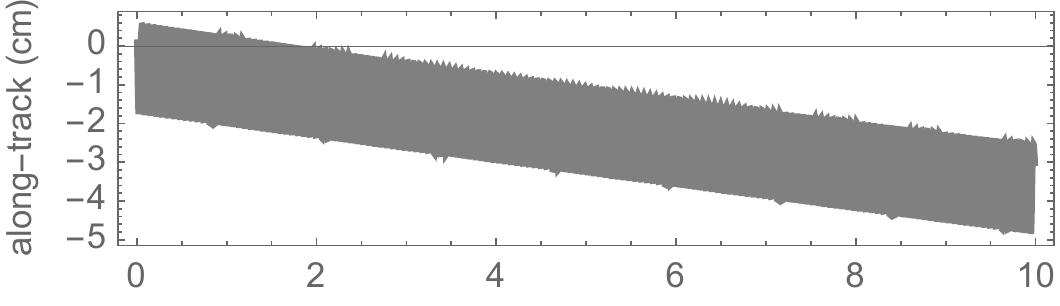}
\includegraphics[scale=0.7]{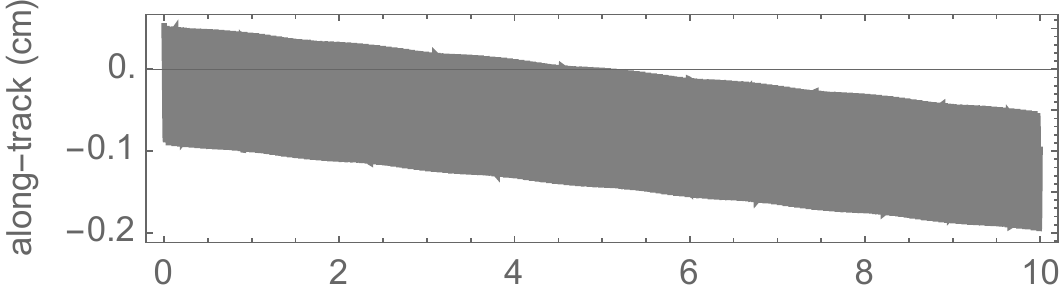}
\caption{Along-track errors of the second order theory of the traditional (top) and extended phase space formulation (bottom plot). Note the different scales of the ordinates. Abscissas are days.}
\label{fig:along232both}
\end{figure}

When the second-order theory is used we obtain notable accuracy gains in both cases, as expected from the characteristics of perturbation solutions. Again, the amplitude of the periodic errors in the radial and cross-track directions is similar for both perturbation theories, now in the mm level, and corresponding errors are not presented. On the contrary, while the amplitude of the periodic errors in the along-track direction remains also in the mm level in the case of the solution in extended phase space, it is not at all the case of the errors of the traditional theory, whose amplitude now reaches the cm level. This is shown in Fig.~\ref{fig:along232both}, where we also notice the secular rate of the along-track errors of the traditional theory of about one third of cm per day, which is much smaller than the rate of about one tenth of mm per day of the extended phase space solution. Nonetheless, these figures yield an analogous ratio to the first-order case.
\par

\begin{figure}[htb]
\centering
\includegraphics[scale=0.7]{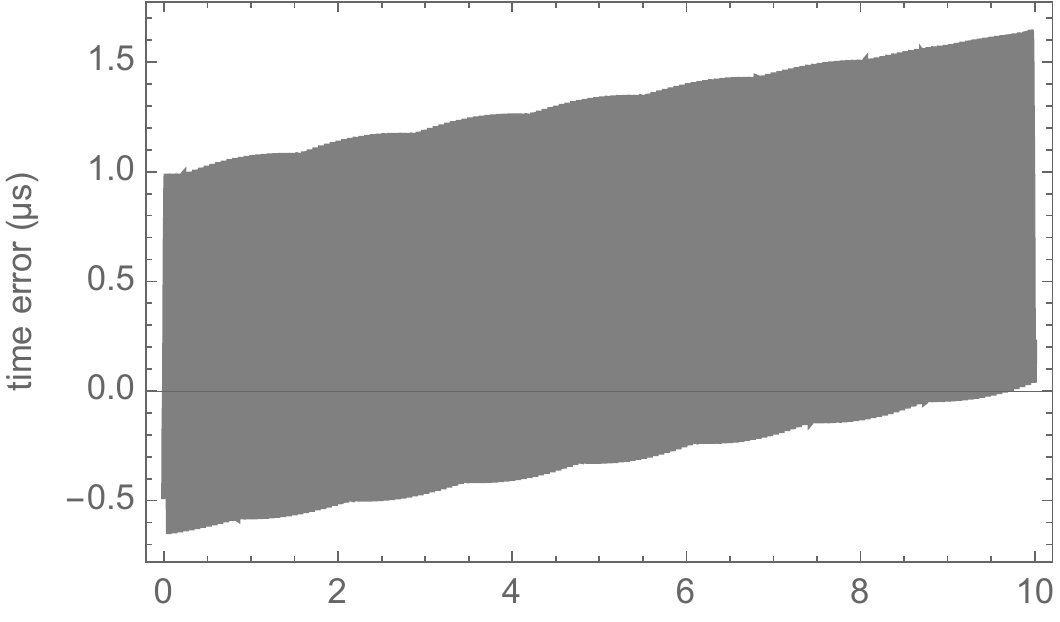}
\caption{Errors in the physical time determination stemming from the second-order theory in the extended phase space. Abscissas are days. }
\label{fig:timerror2}
\end{figure}

\begin{figure}[htb]
\centering
\includegraphics[scale=0.7]{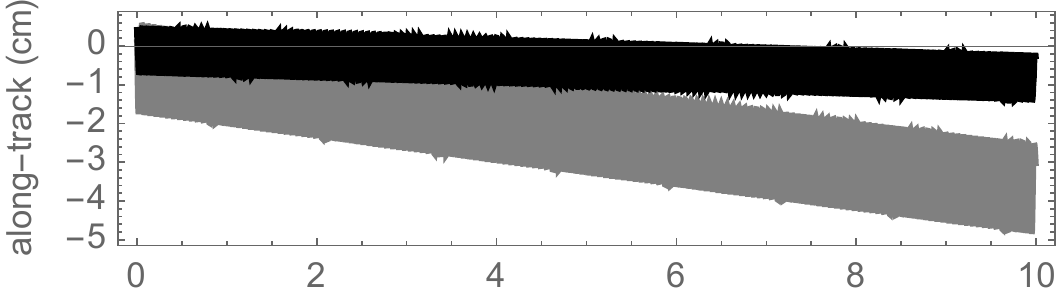}
\caption{Along-track errors of the second-order theory in the extended phase space with the physical time as argument (black) superimposed to the top plot of Fig.~\protect\ref{fig:along232both}. Abscissas are days. }
\label{fig:intrack232time}
\end{figure}

Again, additional inaccuracies to those presented in Fig.~\ref{fig:along232both} arise in the case of the extended phase space solution due to the need of computing the physical time from the fictitious one, or vice versa. As shown in Fig.~\ref{fig:timerror2}, the amplitude of these errors is now reduced to about 1 microsecond, which is not enough to hidden a clear secular trend of about $0.1\,\mathrm{{\mu}s/day}$. These additional errors of the extended phase space formulation, now reaching the cm level for the PRISMA orbit, result in a notably increase with respect to the mm level shown in the bottom plot of Fig.~\ref{fig:along232both}. Again, timing errors notably balance the accuracy performance of both kinds of theories, as illustrated in Fig.~\ref{fig:intrack232time}. Therefore, like in the first-order case, this source of errors may be a dominant part that cannot be ignored in the evaluation performance of perturbation solutions based on the formulation in the extended phase space.
\par

\section{Conclusions}

Inclusion of the total energy among the variables of an analytic orbit generator has beneficial, radical effects in the propagation of along-track errors of circumterrestrial orbits, which are commonly reduced by one order of magnitude with respect to usual along-track errors obtained with traditional (properly initialized) perturbation solutions. This fact could make the extended phase formulation preferable for mission analysis and planning. However, this kind of solution has to deal unavoidably with the transformation between physical and fictitious times in ephemeris computation. This is not only a clear inconvenience from the algorithmic point of view, but, most notably, this change cannot be carried out exactly due to the truncation that is inherent to any perturbation solution. Additional inaccuracies stemming from this conversion, which show both secular and periodic trends, may be as important as those derived from the intrinsic components of the position errors. Therefore, they cannot be ignored in the objective assessment of accuracy performance. From the point of view of the construction of the analytical solution, the fictitious-time approach has clear advantages regarding the computation of the integrals involved in the process, yet a post-processing is needed to reduce the length of the series comprising the solution to a comparable size to those resulting from the traditional approach. These pros and cons do not present a clear advantage for one of the options, and make the choice between traditional and extended phase space orbital perturbation solutions mostly a matter of the particular needs of prospective users.

\section*{Acknowledgments}
Partial support by the Spanish State Research Agency and the European Regional Development Fund (Project PID2020-112576GB-C22, AEI/ERDF, EU) is recognized.


\begin{thebibliography}{10}
\small

\bibitem{Ahmed1994}
M.~K.~M. {Ahmed}.
\newblock {On the normalization of perturbed Keplerian systems}.
\newblock {\em The Astronomical Journal}, 107:1900--1903, May 1994.

\bibitem{AlfriendCoffey1984}
K.~T. {Alfriend} and S.~L. {Coffey}.
\newblock {Elimination of the perigee in the satellite problem}.
\newblock {\em Celestial Mechanics}, 32(2):163--172, February 1984.

\bibitem{BoccalettiPucacco1998v2}
D.~{Boccaletti} and G.~{Pucacco}.
\newblock {\em {Theory of orbits. Volume 2: Perturbative and geometrical
  methods}}.
\newblock Astronomy and Astrophysics Library. Springer-Verlag, Berlin
  Heidelberg New York, 1st edition, 2002.

\bibitem{BonavitoWatsonWalden1969}
N.~L. {Bonavito}, S.~{Watson}, and H.~{Walden}.
\newblock {An Accuracy and Speed Comparison of the Vinti and Brouwer Orbit
  Prediction Methods}.
\newblock Technical Report NASA TN D-5203, Goddard Space Flight Center,
  Greenbelt, Maryland, May 1969.

\bibitem{Bond1979}
V.~R. Bond.
\newblock {An Analytical Singularity-Free Solution to the $J_2$ Perturbation
  Problem}.
\newblock Technical Report NASA-TM-58221; JSC-13128, NASA Johnson Space Center,
  jul 1979.

\bibitem{BreakwellVagners1970}
J.~V. {Breakwell} and J.~{Vagners}.
\newblock {On Error Bounds and Initialization in Satellite Orbit Theories}.
\newblock {\em Celestial Mechanics}, 2:253--264, June 1970.

\bibitem{Brouwer1959}
D.~{Brouwer}.
\newblock {Solution of the problem of artificial satellite theory without
  drag}.
\newblock {\em The Astronomical Journal}, 64:378--397, November 1959.

\bibitem{Cain1962}
B.~J. {Cain}.
\newblock {Determination of mean elements for Brouwer's satellite theory}.
\newblock {\em Astronomical Journal}, 67:391--392, August 1962.

\bibitem{CoffeyAlfriend1984}
S.~Coffey and K.~T. Alfriend.
\newblock {An analytical orbit prediction program generator}.
\newblock {\em Journal of Guidance, Control and Dynamics}, 7(5):575--581, 1984.

\bibitem{CoffeyNealSegermanTravisano1995}
S.~L. Coffey, H.~L. Neal, A.~M. Segerman, and J.~J. Travisano.
\newblock An analytic orbit propagation program for satellite catalog
  maintenance.
\newblock In K.~Terry Alfriend, I.~Michael Ross, Arun~K. Misra, and C.~Fred
  Peters, editors, {\em AAS/AIAA Astrodynamics Conference 1995}, volume~90 of
  {\em Advances in the Astronautical Sciences}, pages 1869--1892, P.O.~Box
  28130, San Diego, California 92198, USA, 1996. American Astronautical
  Society, Univelt, Inc.

\bibitem{Deprit1969}
A.~{Deprit}.
\newblock Canonical transformations depending on a small parameter.
\newblock {\em Celestial Mechanics}, 1(1):12--30, 1969.

\bibitem{Deprit1981}
A.~Deprit.
\newblock The elimination of the parallax in satellite theory.
\newblock {\em Celestial Mechanics}, 24(2):111--153, 1981.

\bibitem{Deprit1982}
A.~{Deprit}.
\newblock {Delaunay normalisations}.
\newblock {\em Celestial Mechanics}, 26:9--21, January 1982.

\bibitem{DepritFerrer1987}
A.~Deprit and S.~Ferrer.
\newblock {Note on Cid's Radial Intermediary and the Method of Averaging}.
\newblock {\em Celestial Mechanics}, 40(3-4):335--343, 1987.

\bibitem{DepritRom1970}
A.~Deprit and A.~Rom.
\newblock {The Main Problem of Artificial Satellite Theory for Small and
  Moderate Eccentricities}.
\newblock {\em Celestial Mechanics}, 2(2):166--206, June 1970.

\bibitem{EcksteinHechler1970}
M.~C. Eckstein and F.~Hechler.
\newblock {A reliable derivation of the perturbations due to any zonal and
  tesseral harmonics of the geopotential for nearly-circular satellite orbits}.
\newblock Scientific Report ESRO SR-13, European Space Research Organisation,
  Darmstadt, Federal Republic of Germany, June 1970.

\bibitem{Floria1997}
Luis {Flor{\'\i}a}.
\newblock {Perturbed Gylden Systems and Time-Dependent Delaunay-Like
  Transformations}.
\newblock {\em Celestial Mechanics and Dynamical Astronomy}, 68(1):75--85, May
  1997.

\bibitem{GaiasColomboLara2020}
G.~{Gaias}, C.~{Colombo}, and M.~{Lara}.
\newblock {Analytical Framework for Precise Relative Motion in Low Earth
  Orbits}.
\newblock {\em Journal of Guidance Control Dynamics}, 43(5):915--927, March
  2020.

\bibitem{Healy2000}
L.~M. {Healy}.
\newblock {The Main Problem in Satellite Theory Revisited}.
\newblock {\em Celestial Mechanics and Dynamical Astronomy}, 76(2):79--120,
  2000.

\bibitem{Jefferys1971}
W.~H. {Jefferys}.
\newblock {Automated, Closed Form Integration of Formulas in Elliptic Motion.}
\newblock {\em Celestial Mechanics}, 3:390--394, September 1971.

\bibitem{Kozai1959}
Y.~{Kozai}.
\newblock {The motion of a close earth satellite}.
\newblock {\em The Astronomical Journal}, 64:367--377, November 1959.

\bibitem{Kozai1962}
Y.~{Kozai}.
\newblock {Second-order solution of artificial satellite theory without air
  drag}.
\newblock {\em The Astronomical Journal}, 67:446--461, September 1962.

\bibitem{Lara2019UR}
M.~{Lara}.
\newblock {Review of analytical solutions for low earth orbit propagation and
  study of the precision improvement in the conversion of osculating to mean
  elements}.
\newblock Technical Report CM 2019/SER0023, Universidad de La Rioja, Logro\~no,
  La Rioja, September 2019.

\bibitem{Lara2020arxiv}
M.~{Lara}.
\newblock {Brouwer's satellite solution redux}.
\newblock {\em Celestial Mechanics and Dynamical Astronomy}, pages 1--20,
  September 2021.

\bibitem{Lara2021IAC}
M.~{Lara}.
\newblock {Improving efficiency of analytic orbit propagation
  (IAC-21,C1,7,2,x65390)}.
\newblock In {\em Proceedings of the 72nd International Astronautical Congress
  (IAC), Dubai, United Arab Emirates, 25--29 October 2021}. International
  Astronautical Federation (IAF), International Astronautical Federation (IAF),
  2021.

\bibitem{LaraSanJuanLopezOchoa2013c}
M.~{Lara}, J.~F. {San-Juan}, and L.~M. {L{\'o}pez-Ochoa}.
\newblock {Delaunay variables approach to the elimination of the perigee in
  Artificial Satellite Theory}.
\newblock {\em Celestial Mechanics and Dynamical Astronomy}, 120(1):39--56,
  September 2014.

\bibitem{LaraSanJuanLopezOchoa2013b}
M.~{Lara}, J.~F. {San-Juan}, and L.~M. {L{\'o}pez-Ochoa}.
\newblock {Proper Averaging Via Parallax Elimination (AAS 13-722)}.
\newblock In Stephen~B. Broschart, James~D. Turner, Kathleen~C. Howell, and
  Felix~R. Hoots, editors, {\em Astrodynamics 2013}, volume 150 of {\em
  {Advances in the Astronautical Sciences}}, pages 315--331, P.O.~Box 28130,
  San Diego, California 92198, USA, January 2014. American Astronautical
  Society, Univelt, Inc.

\bibitem{Lara2019CMDA}
Martin {Lara}.
\newblock {A new radial, natural, higher-order intermediary of the main problem
  four decades after the elimination of the parallax}.
\newblock {\em Celestial Mechanics and Dynamical Astronomy}, 131(9):1--20,
  September 2019.

\bibitem{Lara2020}
Martin {Lara}.
\newblock {Solution to the main problem of the artificial satellite by reverse
  normalization}.
\newblock {\em Nonlinear Dynamics}, 101(2):1501--1524, July 2020.

\bibitem{Lara2021}
Martin Lara.
\newblock {\em {Hamiltonian Perturbation Solutions for Spacecraft Orbit
  Prediction. The method of Lie Transforms}}, volume~54 of {\em De Gruyter
  Studies in Mathematical Physics}.
\newblock De Gruyter, Berlin/Boston, 1 edition, 2021.

\bibitem{Lyddane1963}
R.~H. {Lyddane}.
\newblock {Small eccentricities or inclinations in the Brouwer theory of the
  artificial satellite}.
\newblock {\em {Astronomical Journal}}, 68(8):555--558, October 1963.

\bibitem{LyddaneCohen1962}
R.~H. {Lyddane} and C.~J. {Cohen}.
\newblock {Numerical comparison between Brouwer's theory and solution by
  Cowell's method for the orbit of an artificial satellite}.
\newblock {\em Astronomical Journal}, 67:176--177, April 1962.

\bibitem{Metris1991}
G.~{Metris}.
\newblock {Mean values of particular functions in the elliptic motion}.
\newblock {\em Celestial Mechanics and Dynamical Astronomy}, 52:79--84, March
  1991.

\bibitem{PerssonJacobssonGill2005}
S.~{Persson}, B.~{Jacobsson}, and E.~{Gill}.
\newblock {PRISMA -- Demonstration Mission for Advanced Rendezvous and
  Formation Flying Technologies and Sensors (paper IAC-05-B56B07)}.
\newblock In {\em Proceedings of the 56th International Astronautical Congress
  (IAC), October 17 - 21 2005, Fukuoka, Japan}, pages 1--10. International
  Astronautical Federation (IAF), International Astronautical Federation (IAF),
  October 2005.

\bibitem{Poincare1892vII}
Henri {Poincar\'e}.
\newblock {\em {Les m\'ethodes nouvelles de la m\'ecanique c\'eleste. Tome 2}}.
\newblock Gauthier-Villars et fils (Paris), 1893.

\bibitem{SanJuan1994}
J.~F. San-Juan.
\newblock {ATESAT: Automatization of theories and ephemeris in the artificial
  satellite problem}.
\newblock Technical Report CT/TI/MS/MN/94-250, Centre National d'\'Etudes
  Spatiales, 18, avenue Edouard Belin - 31401 Toulouse Cedex 9, France, May
  1994.

\bibitem{Scheifele1970}
G.~{Scheifele}.
\newblock {G{\'e}n{\'e}ralisation des {\'e}l{\'e}ments de Delaunay en
  m{\'e}canique c{\'e}leste. Application au mouvement d'un satellite
  artificiel.}
\newblock {\em Academie des Sciences Paris Comptes Rendus Serie B Sciences
  Physiques}, 271:729--732, January 1970.

\bibitem{Scheifele1970CeMec}
G.~{Scheifele}.
\newblock {On Nonclassical Canonical Systems}.
\newblock {\em Celestial Mechanics}, 2(3):296--310, September 1970.

\bibitem{Scheifele1981}
G.~{Scheifele}.
\newblock {An analytical singularity-free orbit predictor for near-earth
  satellites.}
\newblock In {\em Proceedings of the International Symposium on Spacecraft
  Flight Dynamics, Darmstadt, 18--22 May 1981}, volume 160 of {\em ESA Special
  Publication}, pages 299--305, August 1981.

\bibitem{ScheifeleGraf1974}
G~{Scheifele} and O~{Graf}.
\newblock {Analytical satellite theories based on a new set of canonical
  elements}.
\newblock In {\em Mechanics and Control of Flight Conference}, pages 1--20,
  Reston, Virigina, feb 1974. American Institute of Aeronautics and
  Astronautics.

\bibitem{StiefelScheifele1971}
E.~L. {Stiefel} and G.~{Scheifele}.
\newblock {\em {Linear and Regular Celestial Mechanics}}, volume 174 of {\em
  Grundlehren der mathematischen Wissenschaften}.
\newblock Springer-Verlag, Berlin Heidelberg, 1 edition, 1971.

\bibitem{Ustinov1967}
B.~A. {Ustinov}.
\newblock {Motion of Satellites in Small-Eccentricity Orbits in the Noncentral
  Gravitational Field of the Earth}.
\newblock {\em Cosmic Research}, 5:159--168, March 1967.

\bibitem{Walter1967}
H.~G. {Walter}.
\newblock {Conversion of osculating orbital elements into mean elements}.
\newblock {\em Astronomical Journal}, 72:994--997, October 1967.

\end{thebibliography}
\end{document}